# Certain Classes of Univalent Functions With Negative Coefficients Defined By General Linear Operator

**Mazin Sh.Mahmoud[1] , Abdul Rahman S.Juma[2], Raheam A. Mansor Al-Saphory[3]**

*[2] Department of Mathematics, College of Education for pure Sciences, University of Tikrit , Tikrit , Iraq*

**Corresponding Author:**

**Name:** Mazin Sh. Mahmoud

**E-mail:** maz2004in@yahoo.com

**Tel:**

## ABSTRACT

In this study, a subclass $S_m^{s,c}(\mu,\beta,\delta)$of an univalent function with negative coefficients which is defined by anew general Linear operator $\mathcal{H}_m^{S,C}$ have been introduced. The sharp results for coefficients estimators, distortion and closure bounds, Hadamard product, and Neighborhood, and this paper deals with the utilizing of many of the results for classical hypergeometric function, where there can be generalized to m-hypergeometric functions.

A subclasses of univalent functions are presented, and it has involving operator $\mathcal{H}_m^{s,c}(c_i,b_j)$which generalizes many well-known. Denote A the class of functions *f* and we have other results have been studied.

## 1. Introduction

Many researchers such as Mohammed and Darus [1], Aldweby and Darus [2],and others have used the m-hypergeometric functions for studying certain families of mathematic viable functions in an open disk unit. The m-hypergeometric functions are generalized configuration of the classical hypergeometric functions. Then by assuming the limit m → 1, it would return to a classical hypergeometric function. The formal set of hypergeometric functions have been used and introduced by many famous researchers were started by Euler in (1748), Gauss (1813) and Cauchy (1852) see (Juma [3]). Also, it was converted a simple notation $\lim_{m\to 1}\frac{1-m^c}{1-m}=c$ into a systematic theory of hypergemetric function in same trend of theory of Gauss hypergeometric function.

Here, this study deals with the utilizing of many of the results for classical hypergeometric function, where there can be generalized to m-hypergeometric functions.

In this work, a subclasses of univalent functions are introduced, and it has involving operator $\mathcal{H}_m^{s,c}(c_i,b_j)$which generalizes many well-known. Denote A the class of functions *f* of the form

$f(z)=z+\sum_{n=1}^{\infty}a_n z^n$ , (1.1)

which are analytic and univalent in the open unit disk Û={z∈ ℂ:|z|<1}.

A function $f$ ∈A is said to be starlike of complex order if the following condition (see[4]) is satisfied:

$$\Re\left\{\frac{\frac{zf(z))'}{f(z)}-1}{2\delta\left(\frac{zf(z))'}{f(z)}-\mu\right)-\left(\frac{zf(z))'}{f(z)}-1\right)}\right\}>\beta,\ (0\le\mu<\frac{1}{2\delta}\ ,0<\beta\le 1,\frac{1}{2}\le\delta\le 1)\quad (1.2)$$

For complex parameters $c_1,\ldots,c_t$ and $b_1,\ldots,b_r$ ( $b_j\in\mathbb{C}\backslash\{0,-1,-2,\ldots\}$,

j =1,…,r,|*m*|<1), the m-hypergeometric

$${}_t\Psi_r=\sum_{n=0}^{\infty}\frac{(c_1,m)_n\ldots(c_t,m)_n}{(m,m)_n\,(b_1,m)_n\ldots(b_r,m)_n}z^n\qquad(1.3)$$

$(t=r+1$ *such that* $t,r\in N_0=\{0,1,2,3,4\ldots\}$ ; Z ∈ Û).

The m-shifted factorial is involving by

$(c,m)_0=1$ and $(c,m)_n=(1-c)(1-cm)(1-cm^2)\ldots(1-cm^{n-1})$, $n\in\mathbb{N}$,

where c any complex number and in terms of the Gamma function

$$(m^{\mu},m)_n=\frac{\Gamma_m(\mu+n)(1-m)^n}{\Gamma_m(\mu)},$$

such that

$$\Gamma_m(y)=\frac{(m,m)_\infty(1-m)(1-y)}{(m^y,m)_\infty},0<m<1.$$

The study suggests that note that and by utilizing ratio test , the series (1.3)converges absolutely in open unit disk Û, |m|<1

$${}_2\Psi_1=\sum_{n=0}^{\infty}\frac{(c_1,m)_n(c_2,m)_n}{(m,m)_n(b_1,m)_n}z^n\ (|m|<1,z\in\hat{U})$$

Is the m-Gauss hypergeometric function see [4],[5].

Recently Mohammed and Darus [1] defined the following:

$I(c_i; b_j; m)f : A \to A$

$I(c_i; b_j; m)f(z)=z+\sum_{n=2}^{\infty}\frac{(c_1,m)_{n-1}\cdots(c_t,m)_{n-1}}{(m,m)_{n-1}(b_1,m)_{n-1}\cdots(b_r,m)_{n-1}}a_n z^n$ .

The Srivastava-Attiya operator $T_{s,c}$: A $\to$ A is defined in [6] as:

$$T_{s,c}f(z) = z + \sum_{n=2}^{\infty}\left(\frac{1+c}{n+c}\right)^s a_n z^n, \quad (1.4)$$

where $z \in \hat{U}$ , $c \in \mathbb{C}/\{0, -1, -2, \ldots\}$, $s \in \mathbb{C}$ and $f \in A$. This linear operator $\mathbf{T_{s,c}}$ can be written as

$T_{s,c}f(z)=G_{s,c}(z) * f(Z) = (1+c)^s(\Phi(z,s,c)-c^{-s})*f(z)$, by utilizing the Hadamard product (convolution).Here,

$\Phi(z,s,c)=\sum_{n=0}^{\infty}\frac{z^n}{(n+c)^s}$ ,

is the well-known Hurwitz –Lerch zeta function (see[6],[7]). It is also an important function of Analytic Number Theory such the De-Jonquiere function:

$H_{is}(Z)=\sum_{n=0}^{\infty}\frac{z^n}{(n)^s} = z\Phi(z,s,1), \ (Re(s) > 1 \ if \ |z| = 1)$.

We can define the linear operator $\mathcal{H}_m^{s,c}(c_i, b_j)(f)$ : A $\to$ A as follows:

$$\mathcal{H}_m^{s,c}(c_i,b_j)f(z) = z + \sum_{n=2}^{\infty}\frac{(c_1,m)_{n-1}\cdots(c_t,m)_{n-1}}{(m,m)_{n-1}(b_1,m)_{n-1}\cdots(b_r,m)_{n-1}}\left(\frac{1+c}{n+c}\right)^s a_n z^n. \quad (1.5)$$

($z \in U$, $c \in \mathbb{C}/\{0, -1, -2, \ldots\}$, $s \in \mathbb{C}$, $c_i$ , $b_j \in \mathbb{C}/\{0,-1,-2,-3,\ldots\}$, $|m|<1$ and $t=r+1$.

It should be noted that the liner operator (1.5)introduced by A. R.S.Juma and M. Darus[3] .

**Definition 1.** $f$ is a function and $f \in \hat{U}$ is said to be in the class $\mathcal{R}_m^{s,c}(\mu,\beta,\delta)$ if the following condition holds:

$$\left|\frac{\frac{z(\mathcal{H}_m^{s,c}(c_i,b_j)f(z))'}{\mathcal{H}_m^{s,c}(c_i,b_j)f(z)}-1}{2\delta\left(\frac{z(\mathcal{H}_m^{s,c}(c_i,b_j)f(z))'}{\mathcal{H}_m^{s,c}(c_i,b_j)f(z)}-\mu\right)-\left(\frac{z(\mathcal{H}_m^{s,c}(c_i,b_j)f(z))'}{\mathcal{H}_m^{s,c}(c_i,b_j)f(z)}-1\right)}\right| < \beta \quad (1.6)$$

where $0 \le \mu < \frac{1}{2\delta}$ , $0 < \beta \le 1$ , $\frac{1}{2} \le \delta \le 1$ , $z \in U$.

Let T denote the subclass of A consisting of function of the form

$$f(z) = z - \sum_{n=2}^{\infty} a_n z^n, \quad a_n \ge 0 . \quad (1.7)$$

Now we define the class $S_m^{s,c}(\mu,\beta,\delta)$ by

$S_m^{s,c}(\mu,\beta,\delta)= \mathcal{R}_m^{s,c}(\mu,\beta,\delta) \cap T$.

The study have the following class and confirms that note that by specializing the parameters $\mu, \beta, \delta$

1. The class $S_m^{-k,0}(\alpha,\beta,\xi)$ is the class studied by A. R.S.Juma and S. R. Kulkarni [8].

2. The class $S_m^{-k,0}(0,1,1)$ is precisely the class of starlike function in $\hat{U}$.

3.The class $S_m^{-k,0}(\mu,1,1)$ is the class of starlike function of order $\mu$ ($0 \le \mu < 1$).

4.The class $S_m^{-k,0}(0,\beta,\frac{\mu+1}{2})$ is the class studied by Lakshminar-simhan[9].

5. The class $S_m^{-k,c}(\mu,\beta,\delta)$ is the class studied by S. R. Kulkarni [10].

## 2. Conflicients estimates and Other properties

**Theorem 1.** Let $f$ be defined by (1.7).Then $f \in S_m^{s,c}(\mu,\beta,\delta)$ if and only if

$$\sum_{n=2}^{\infty}\left(\frac{(c_1,m)_{n-1}\cdots(c_t,m)_{n-1}}{(m,m)_{n-1}(b_1,m)_{n-1}\cdots(b_r,m)_{n-1}}\right)\left(\frac{1+c}{n+c}\right)^s[(n-1)(1-\beta)2\beta\delta(n-\mu)]a_n \le 2\beta\delta(1-\mu) \quad (2.1)$$

$0 < \beta \le 1, 0 \le \mu < {}^1\!/_2\,\delta, \ \frac{1}{2} \le \delta \le 1.$

**Proof**: If $|z|=1$ ,then

$|z(\mathcal{H}_m^{s,c}(c_i,b_j)f(z))' - (\mathcal{H}_m^{s,c}(c_i,b_j)f(z))|$

$-\beta|2\delta(z(\mathcal{H}_m^{s,c}(c_i,b_j)f(z))' - \mu(\mathcal{H}_m^{s,c}(c_i,b_j)f(z))$

$-(z\mathcal{H}_m^{s,c}(c_i,b_j)f(z))' - \mathcal{H}_m^{s,c}(c_i,b_j)f(z))|$.

By utilizing (1.5)we have

$\mathcal{H}_m^{s,c}(c_i,b_j)f(z))' = z + \sum_{n=2}^{\infty} n\left(\frac{(c_1,m)_{n-1}\cdots(c_t,m)_{n-1}}{(m,m)_{n-1}(b_1,m)_{n-1}\cdots(b_r,m)_{n-1}}\right)\left(\frac{1+c}{n+c}\right)^s a_n z^{n-1}$

$z\mathcal{H}_m^{s,c}(c_i,b_j)f(z))' = z + \sum_{n=2}^{\infty} n\left(\frac{(c_1,m)_{n-1}\cdots(c_t,m)_{n-1}}{(m,m)_{n-1}(b_1,m)_{n-1}\cdots(b_r,m)_{n-1}}\right)\left(\frac{1+c}{n+c}\right)^s a_n z^{n}$

$= |z + \sum_{n=2}^{\infty} n\left(\frac{(c_1,m)_{n-1}\cdots(c_t,m)_{n-1}}{(m,m)_{n-1}(b_1,m)_{n-1}\cdots(b_r,m)_{n-1}}\right)\left(\frac{1+c}{n+c}\right)^s a_n z^{n}|$

$-\beta(2\delta z + \sum_{n=2}^{\infty} n\left(\frac{(c_1,m)_{n-1}\cdots(c_t,m)_{n-1}}{(m,m)_{n-1}(b_1,m)_{n-1}\cdots(b_r,m)_{n-1}}\right)\left(\frac{1+c}{n+c}\right)^s a_n z^{n})$

$-\mu z - \sum_{n=2}^{\infty}\left(\frac{(c_1,m)_{n-1}\cdots(c_t,m)_{n-1}}{(m,m)_{n-1}(b_1,m)_{n-1}\cdots(b_r,m)_{n-1}}\right)\left(\frac{1+c}{n+c}\right)^s a_n z^{n}$

$-z-\sum_{n=2}^{\infty} n\left(\frac{(c_1,m)_{n-1}\cdots(c_t,m)_{n-1}}{(m,m)_{n-1}(b_1,m)_{n-1}\cdots(b_r,m)_{n-1}}\right)\left(\frac{1+c}{n+c}\right)^s a_n z^{n} + z$

$+ \sum_{n=z}^{\infty}\left(\frac{(c_1,m)_{n-1}\cdots(c_t,m)_{n-1}}{(m,m)_{n-1}(b_1,m)_{n-1}\cdots(b_r,m)_{n-1}}\right)\left(\frac{1+c}{n+c}\right)^s a_n z^{n}|$

$=|\sum_{n=2}^{\infty}(n-1)\left(\frac{(c_1,m)_{n-1}\cdots(c_t,m)_{n-1}}{(m,m)_{n-1}(b_1,m)_{n-1}\cdots(b_r,m)_{n-1}}\right)\left(\frac{1+c}{n+c}\right)^s a_n z^{n}| - \beta|2\delta z + \sum_{n=2}^{\infty} 2\delta n\left(\frac{(c_1,m)_{n-1}\cdots(c_t,m)_{n-1}}{(m,m)_{n-1}(b_1,m)_{n-1}\cdots(b_r,m)_{n-1}}\right)\left(\frac{1+c}{n+c}\right)^s a_n z^{n}$

$- \mu\delta Z - \sum_{n=2}^{\infty} 2\delta\mu\left(\frac{(c_1,m)_{n-1}\cdots(c_t,m)_{n-1}}{(m,m)_{n-1}(b_1,m)_{n-1}\cdots(b_r,m)_{n-1}}\right)\left(\frac{1+c}{n+c}\right)^s a_n z^{n}$

$- \sum_{n=2}^{\infty} n\left(\frac{(c_1,m)_{n-1}\cdots(c_t,m)_{n-1}}{(m,m)_{n-1}(b_1,m)_{n-1}\cdots(b_r,m)_{n-1}}\right)\left(\frac{1+c}{n+c}\right)^s a_n z^{n} + \sum_{n=2}^{\infty}\left(\frac{(c_1,m)_{n-1}\cdots(c_t,m)_{n-1}}{(m,m)_{n-1}(b_1,m)_{n-1}\cdots(b_r,m)_{n-1}}\right)\left(\frac{1+c}{n+c}\right)^s a_n z^{n}|$

$=|\sum_{n=2}^{\infty}(n-1)\left(\frac{(c_1,m)_{n-1}\cdots(c_t,m)_{n-1}}{(m,m)_{n-1}(b_1,m)_{n-1}\cdots(b_r,m)_{n-1}}\right)\left(\frac{1+c}{n+c}\right)^s a_n z^{n}$

$-\beta|2\delta z(1-\mu) + 2\delta\sum_{n=2}^{\infty}(n-\mu)\left(\frac{(c_1,m)_{n-1}\cdots(c_t,m)_{n-1}}{(m,m)_{n-1}(b_1,m)_{n-1}\cdots(b_r,m)_{n-1}}\right)\left(\frac{1+c}{n+c}\right)^s a_n z^{n}$

$- \sum_{n=2}^{\infty}(n-\mu)\left(\frac{(c_1,m)_{n-1}\cdots(c_t,m)_{n-1}}{(m,m)_{n-1}(b_1,m)_{n-1}\cdots(b_r,m)_{n-1}}\right)\left(\frac{1+c}{n+c}\right)^s a_n z^{n}$

$\le \sum_{n=2}^{\infty}\left(\frac{(c_1,m)_{n-1}\cdots(c_t,m)_{n-1}}{(m,m)_{n-1}(b_1,m)_{n-1}\cdots(b_r,m)_{n-1}}\right)\left(\frac{1+c}{n+c}\right)^s[(n-1)(1-\beta) + 2\beta\delta(n-\mu)]a_n - 2\beta\delta(1-\mu)] \le 0$

$\sum_{n=2}^{\infty}\left(\frac{(c_1,m)_{n-1}\cdots(c_t,m)_{n-1}}{(m,m)_{n-1}(b_1,m)_{n-1}\cdots(b_r,m)_{n-1}}\right)\left(\frac{1+c}{n+c}\right)^s[(n-1)(1-\beta) + 2\beta\delta(n-\mu)]a_n \le 2\beta\delta(1-\mu)$

By hypothesis thus by maximum modulus theorem , we get $f \in S_m^{s,c}(\mu,\beta,\delta)$.

And versa, suppose that $f\in S_m^{s,c}(\mu,\beta,\delta)$,therefore the condition (1.7) gives us

$$\left|\frac{\frac{z(\mathcal{H}_m^{s,c}(c_i,b_j)f(z))'}{\mathcal{H}_m^{s,c}(c_i,b_j)f(z)}-1}{2\delta\left(\frac{z(\mathcal{H}_m^{s,c}(c_i,b_j)f(z))'}{\mathcal{H}_m^{s,c}(c_i,b_j)f(z)}-\mu\right)-\left(\frac{z(\mathcal{H}_m^{s,c}(c_i,b_j)f(z))'}{\mathcal{H}_m^{s,c}(c_i,b_j)f(z)}-1\right)}\right|<\beta$$

$$=\left|\left[-\sum_{n=2}^{\infty}\left[\left(\frac{(c_1,m)_{n-1}\dots(c_t,m)_{n-1}}{(m,m)_{n-1}(b_1,m)_{n-1}\dots(b_r,m)_{n-1}}\right)\left(\frac{1+c}{n+c}\right)^s\right](n-1)a_nz^{n-1}\right]/[(2\delta(1-\mu)-2\delta)\times 2\sum_{n=2}^{\infty}(n-\mu)\left(\frac{(c_1,m)_{n-1}\dots(c_t,m)_{n-1}}{(m,m)_{n-1}(b_1,m)_{n-1}\dots(b_r,m)_{n-1}}\right)\left(\frac{1+c}{n+c}\right)^s a_nz^{n-1}+\sum_{n=2}^{\infty}\left(\frac{(c_1,m)_{n-1}\dots(c_t,m)_{n-1}}{(m,m)_{n-1}(b_1,m)_{n-1}\dots(b_r,m)_{n-1}}\right)\left(\frac{1+c}{n+c}\right)^s(n-1)a_nz^{n-1}\right|<\beta.$$

Sine |Re(z)|<|z| for all z , we have

$$\mathrm{Re}\left\{\left[\sum_{n=2}^{\infty}\left(\frac{(c_1,m)_{n-1}\dots(c_t,m)_{n-1}}{(m,m)_{n-1}(b_1,m)_{n-1}\dots(b_r,m)_{n-1}}\right)\left(\frac{1+c}{n+c}\right)^s(n-1)a_nz^{n-1}\right]/[(2\delta(1-\mu)-2\delta)\times\left(\sum_{n=2}^{\infty}(n-\mu)\left(\frac{(c_1,m)_{n-1}\dots(c_t,m)_{n-1}}{(m,m)_{n-1}(b_1,m)_{n-1}\dots(b_r,m)_{n-1}}\right)\left(\frac{1+c}{n+c}\right)^s a_nz^{n-1}\right)+\sum_{n=2}^{\infty}\left(\frac{(c_1,m)_{n-1}\dots(c_t,m)_{n-1}}{(m,m)_{n-1}(b_1,m)_{n-1}\dots(b_r,m)_{n-1}}\right)\left(\frac{1+c}{n+c}\right)^s(n-1)a_nz^{n-1}]\right\}<\beta.$$

Let $z\to1^{-}$through real values . Then we get (2.1) the result is sharp for function

$$f(z)=z-\frac{2\beta\delta(1-\mu)}{(n-1)(-\beta+1)+2\beta\delta(n-\mu)[\left(\frac{(c_1,m)_{n-1}\dots(c_t,m)_{n-1}}{(m,m)_{n-1}(b_1,m)_{n-1}\dots(b_r,m)_{n-1}}\right)\left(\frac{1+c}{n+c}\right)^s]}z^n,n\geq 2.$$

□

**Corollary2.1:** Let $f$ belong to the class $S_m^{s,c}(\mu,\beta,\delta)$.Then

$$a_n\leq\frac{2\beta\delta(1-\mu)}{(1-\beta)(n-1)+2\beta\delta(n-\mu)\left(\frac{(c_1,m)_{n-1}\dots(c_t,m)_{n-1}}{(m,m)_{n-1}(b_1,m)_{n-1}\dots(b_r,m)_{n-1}}\right)\left(\frac{1+c}{n+c}\right)^s},n\geq 2.\quad(2.2)$$

□

**Theorem 2.** Let $f\in S_m^{s,c}(\mu,\beta,\delta)$. Then for |z|≤ r < 1 , we get

$$r-\frac{r^2(2\beta\delta(1-\mu))}{\left(\frac{(c_1,m)\dots(c_t,m)}{(m,m)(b_1,m)\dots(b_r,m)}\right)\left(\frac{1+c}{2+c}\right)^s[(1-\beta)+2\beta\delta(2-\mu)]}\leq\left|\mathcal{H}_m^{s,c}(c_i,b_j)Hz\right|\leq r+r^2\frac{[(2\beta\delta(1-\mu)]}{\left(\frac{(c_1,m)\dots(c_t,m)}{(m,m)(b_1,m)\dots(b_r,m)}\right)\left(\frac{1+c}{2+c}\right)^s[(1-\beta)+2\beta\delta(2-\mu)]}\quad(2.3)$$

$$1-2r\frac{(2\beta\delta(1-\mu))}{\left(\frac{(c_1,m)\dots(c_t,m)}{(m,m)(b_1,m)\dots(b_r,m)}\right)\left(\frac{1+c}{2+c}\right)^s[2\beta\delta(2-\mu)+(1-\beta)]}\leq\left|\mathcal{H}_m^{S,C}(c_i,b_j)f(z))'\right|\leq 1+\frac{(2\beta\delta(1-\mu))}{\left(\frac{(c_1,m)\dots(c_t,m)}{(m,m)(b_1,m)\dots(b_r,m)}\right)\left(\frac{1+c}{2+c}\right)^s[2\beta\delta(2-\mu)+(1-\beta)]}\quad(2.4)$$

The above bounds are sharp.

**Proof.** By theorem 1 , we have

$$\sum_{n=2}^{\infty}\left(\frac{(c_1,m)_{n-1}\dots(c_t,m)_{n-1}}{(m,m)_{n-1}(b_1,m)_{n-1}\dots(b_r,m)_{n-1}}\right)\left(\frac{1+c}{n+c}\right)^s[(1-\beta)(n-1)-2\beta\delta(\mu-n)]a_n\leq 2\beta\delta(1-\mu),$$

then we have

$$\left(\frac{(c_1,m)\dots(c_t,m)}{(m,m)(b_1,m)\dots(b_r,m)}\right)\left(\frac{1+c}{2+c}\right)^s[(1-\beta)+2\beta\delta(2-\mu)]a_n\leq\sum_{n=2}^{\infty}\left(\frac{(c_1,m)_{n-1}\dots(c_t,m)_{n-1}}{(m,m)_{n-1}(b_1,m)_{n-1}\dots(b_r,m)_{n-1}}\right)\left(\frac{1+c}{n+c}\right)^s[(1-\beta)(n-1)-2\beta\delta(\mu-n)]a_n\leq 2\beta\delta(1-\mu).$$

Then

$$\sum_{n=2}^{\infty}a_n\leq\frac{2\beta\delta(1-\mu)}{\left(\frac{(c_1,m)\dots(c_t,m)}{(m,m)(b_1,m)\dots(b_r,m)}\right)\left(\frac{1+c}{2+c}\right)^s[(1-\beta)+2\beta\delta(2-\mu)]}.$$

Hence

$$\left|\mathcal{H}_m^{S,C}(c_i,b_j)f(z))'\right|\leq |z|+|z|^2\left(\frac{(c_1,m)\dots(c_t,m)}{(m,m)(b_1,m)\dots(b_r,m)}\right)\left(\frac{1+c}{2+c}\right)^s\sum_{n=2}^{\infty}a_n\leq r+r^2\left[\left(\frac{(c_1,m)\dots(c_t,m)}{(m,m)(b_1,m)\dots(b_r,m)}\right)\right]\left(\frac{1+c}{2+c}\right)^s\sum_{n=2}^{\infty}a_n\leq r+\frac{2r^2\beta\delta(1-\mu)}{(1-\beta)+2\beta\delta(2-\mu)},$$

and

$$\left|\mathcal{H}_m^{S,C}(c_i,b_j)f(z))'\right|\geq r-r^2\left[\frac{(c_1,m)\dots(c_t,m)}{(m,m)(b_1,m)\dots(b_r,m)}\left(\frac{1+c}{2+c}\right)^s\right]\sum_{n=2}^{\infty}a_n\geq r-\frac{2r^2\beta\delta(1-\mu)}{(1-\beta)+2\beta\delta(2-\mu)},$$

thus (2.3)is true. Further

$$\left|\mathcal{H}_m^{S,C}(c_i,b_j)f(z))'\right|\leq 1+2r\left[\frac{(c_1,m)\dots(c_t,m)}{(m,m)(b_1,m)\dots(b_r,m)}\left(\frac{1+c}{2+c}\right)^s\right]\sum_{n=2}^{\infty}a_n\leq 1+\frac{4r\beta\delta(1-\mu)}{2\beta\delta(2-\mu)+(1-\beta)}.$$

And also

$$\left|\mathcal{H}_m^{S,C}(c_i,b_j)f(z))'\right|\geq 1-\frac{4r\beta\delta(1-\mu)}{2\beta\delta(2-\mu)+(1-\beta)}.$$

The result is sharp for function $f(z)$,defined by

$$f(z)=z+\frac{2\beta\delta(1-\mu)}{(1-\beta)+[2\beta\delta(2-\mu)]}z^2,\ z=\pm r.$$

□

**Theorem3. Let** $0<\beta\leq 1, 0<\mu_1\leq\mu_2<\frac{1}{2\delta}$ *and* $\frac{1}{2}\leq\delta\leq 1$ *the* $S_m^{s,c}(\mu_2,\beta,\delta)\subset S_m^{s,c}(\mu_1,\beta,\delta)$.

**Proof:** By utilizing assumption we get

$$\frac{2\beta\delta(1-\mu_2)}{\left(\frac{(c_1,m)\dots(c_t,m)}{(m,m)(b_1,m)\dots(b_r,m)}\right)\left(\frac{1+c}{2+c}\right)^s[(n-1)(1-\beta)-2\beta\delta(\mu_2-n)]}\leq\frac{2\beta\delta(1-\mu_1)}{\left(\frac{(c_1,m)\dots(c_t,m)}{(m,m)(b_1,m)\dots(b_r,m)}\right)\left(\frac{1+c}{2+c}\right)^s[(n-1)(1-\beta)-2\beta\delta(\mu_1-n)]}.$$

Thus ,$f\in S_m^{s,c}(\mu_1,\beta,\delta)$ implies that

$$\sum_{n=2}^{\infty}\left(\frac{(c_1,m)_{n-1}\dots(c_t,m)_{n-1}}{(m,m)_{n-1}(b_1,m)_{n-1}\dots(b_r,m)_{n-1}}\right)\left(\frac{1+c}{n+c}\right)^s a_n\leq\frac{2\beta\delta(1-\mu_2)}{(1-\beta)(n-1)+2\beta\delta(n-\mu_2)}\leq\frac{2\beta\delta(1-\mu_1)}{(1-\beta)(n-1)+2\beta\delta(n-\mu_1)}$$

then $f\in S_m^{s,c}(\mu_1,\beta,\delta)$

□

**Theorem4**. The set $S_m^{s,c}(\mu,\beta,\delta)$ is the convex set.

**Proof.** Let $f_i(z)=z+\sum_{n=2}^{\infty}a_{n,i}z^n$ (i=1,2) belong to $S_m^{s,c}(\mu,\beta,\delta)$ and let

$g(Z)=\delta_1F_1(Z)+\delta_2F_2(Z)$ with $\delta_1$ and $\delta_2$ no negative and $\delta_1+\delta_2=1$and we write

$$g(z)=z-\sum_{n=2}^{\infty}(\delta_1,a_{n,1}+\delta_2a_{n,2})z^n\ .$$

It is sufficient to show that $g(z)\in S_m^{s,c}(\mu,\beta,\delta)$ that mean

$\sum_{n=2}^{\infty}(\frac{(c_1,m)_{n-1}....(c_t,m)_{n-1}}{(m,m)_{n-1}(b_1,m)_{n-1}...(b_r,m)_{n-1}})(\frac{1+c}{n+c})^s[(1-\beta)(n-1)+2\beta\delta(n-\mu)][\delta_1, a_{n,1}+\delta_2 a_{n,2}]$

$=\delta_1\sum_{n=2}^{\infty}(\frac{(c_1,m)_{n-1}....(c_t,m)_{n-1}}{(m,m)_{n-1}(b_1,m)_{n-1}...(b_r,m)_{n-1}})(\frac{1+c}{n+c})^s[(1-\beta)(n-1)+2\beta\delta(n-\mu)][a_{n,1}]$

$+\delta_2\sum_{n=2}^{\infty}(\frac{(c_1,m)_{n-1}....(c_t,m)_{n-1}}{(m,m)_{n-1}(b_1,m)_{n-1}...(b_r,m)_{n-1}})(\frac{1+c}{n+c})^s[(1-\beta)(n-1)+2\beta\delta(n-\mu)][a_{n,2}]$

$\leq\delta_2(2\beta\delta(1-\mu)+\delta_1(2\beta\delta(1-\mu))=(\delta_2+\delta_2)(2\beta\delta(1-\mu)=2\beta\delta(1-\mu).$

Thus $g(z))\in S_m^{s,c}(\mu,\beta,\delta)$.

The study shall further try to obtain the extreme points in the following theorem .

□

**Theorem5.** Let $f_1$(z) =z and

$f_n(z)=z+\frac{2\beta\delta(1-\mu)}{(\frac{(c_1,m)_{n-1}....(c_t,m)_{n-1}}{(m,m)_{n-1}(b_1,m)_{n-1}...(b_r,m)_{n-1}})(\frac{1+c}{n+c})^s[(1-\beta)(n-1)+2\beta\delta(n-\mu)}z^n$ ,

for all n=2,3,….; $0<\beta\leq 1$, $0\leq\mu<\frac{1}{2\delta}$ , $\frac{1}{2}\leq\delta\leq 1$.

Then $f(z)$ is in the class $S_m^{s,c}(\mu,\beta,\delta)$if and only it can be expressed in the from

$f$(z)=$\sum_{n=1}^{\infty}\gamma_n z^n$ where ($\gamma\geq 0$ *and* $\sum_{n=1}^{\infty}\gamma_n=1$ or $1=\gamma_1+\sum_{n=2}^{\infty}\gamma_n$ ).

**Proof** .Let $f$(z)=$\sum_{n=1}^{\infty}\gamma_n z^n$ where ($\gamma_n\geq 0$ and $\sum_{n=1}^{\infty}\gamma_n=1$).

$f(z)=z+\sum_{n=1}^{\infty}\frac{2\beta\delta(1-\mu)}{(\frac{(c_1,m)_{n-1}....(c_t,m)_{n-1}}{(m,m)_{n-1}(b_1,m)_{n-1}...(b_r,m)_{n-1}})(\frac{1+c}{n+c})^s[(1-\beta)(n-1)+2\beta\delta(n-\mu)}\gamma_n z^n$

, and we obtain

$\sum_{n=1}^{\infty}\left[\frac{(\frac{(c_1,m)_{n-1}....(c_t,m)_{n-1}}{(m,m)_{n-1}(b_1,m)_{n-1}...(b_r,m)_{n-1}})(\frac{1+c}{n+c})^s[(1-\beta)(n-1)+2\beta\delta(n-\mu)}{2\beta\delta(1-\mu)}\right]\times$

$\gamma_n\frac{2\beta\delta(1-\mu)}{(\frac{(c_1,m)_{n-1}....(c_t,m)_{n-1}}{(m,m)_{n-1}(b_1,m)_{n-1}...(b_r,m)_{n-1}})(\frac{1+c}{n+c})^s[(1-\beta)(n-1)+2\beta\delta(n-\mu)}$

$=\sum_{n=1}^{\infty}\gamma_n=1-\gamma_1\leq 1.$

In view of theorem 1, this shows that $f$(z) $\in S_m^{s,c}(\mu,\beta,\delta)$.

Conversely ,

$a_n\leq\frac{2\beta\delta(1-\mu)}{(\frac{(c_1,m)_{n-1}....(c_t,m)_{n-1}}{(m,m)_{n-1}(b_1,m)_{n-1}...(b_r,m)_{n-1}})(\frac{1+c}{n+c})^s[(1-\beta)(n-1)+2\beta\delta(n-\mu)}$ , $n\geq 2$

if

$\gamma_n=\frac{(\frac{(c_1,m)_{n-1}....(c_t,m)_{n-1}}{(m,m)_{n-1}(b_1,m)_{n-1}...(b_r,m)_{n-1}})(\frac{1+c}{n+c})^s[(1-\beta)(n-1)+2\beta\delta(n-\mu)}{2\beta\delta(1-\mu)}$

and $\gamma_1=1-\sum_{n=1}^{\infty}\gamma_n$, then we get

$f(z)=\gamma_1 f_1(z)+\sum_{n=1}^{\infty}\gamma_n f_n(z).$

□

## 3. Neighbourhood and Hadamard product properties

**Definition 3.1[11]:** Let $\gamma\geq 0$ , $f(z)\in T$ on the (1.7 ) the (k, $\gamma$)- neighborhod of a function $f(z)$ defined by

$N_n,\gamma(f)=\{g\in T: g(z)=z-\sum_{n=2}^{\infty}b_n z^n$ and $\sum_{n=2}^{\infty}n\,|a_n-b_n|\leq\gamma\}$, (3.1)

For the identity function *e* (z)=z , we get

$N_n,\gamma(e)=\{g\in T: g(z)=z-\sum_{n=2}^{\infty}b_n z^n$ and $\sum_{n=2}^{\infty}n\,|b_n|\leq\gamma\}$. (3.2)

**Theorem 6.** Let

$\gamma=\frac{4\beta\delta(1-\mu)}{[\frac{(c_1,m)....(c_t,m)}{(m,m)(b_1,m)...(b_r,m)})(\frac{1+c}{2+c})^s](1-\beta)+2\beta\delta(2-\mu)]}$ .

Then $S_m^{s,c}(\mu,\beta,\delta)\subset N_{n,\gamma}(e)$.

**Proof**. Let $f\in S_m^{s,c}(\mu,\beta,\delta)$. Then we get

$((1-\beta)+2\beta\delta(2-\mu))[\frac{(c_1,m)....(c_t,m)}{(m,m)(b_1,m)...(b_r,m)})(\frac{1+c}{2+c})^s]\sum_{n=2}^{\infty}a_n$

$\leq\sum_{n=2}^{\infty}\frac{(c_1,m)_{n-1}....(c_t,m)_{n-1}}{(m,m)_{n-1}(b_1,m)_{n-1}...(b_r,m)_{n-1}})(\frac{1+c}{n+c})^s(1-\beta)(n-1)+2\beta\delta(n-\mu)\leq 2\beta\delta(1-\mu)$,

therefore,

$\sum_{n=2}^{\infty}a_n\leq\frac{2\beta\delta(1-\mu)}{[\frac{(c_1,m)....(c_t,m)}{(m,m)(b_1,m)...(b_r,m)})(\frac{1+c}{2+c})^s](1-\beta)+2\beta\delta(2-\mu)]}$ , (3.3)

also we get $|z|<r$

$|f'(z)|\leq 1+|z|\sum_{n=2}^{\infty}n\,a_n\leq 1+r\sum_{n=2}^{\infty}n\,a_n$.

In view of (3.3)we get

$|f'(z)|\leq 1+r\frac{2(2\beta\delta(1-\mu))}{[\frac{(c_1,m)....(c_t,m)}{(m,m)(b_1,m)...(b_r,m)})(\frac{1+c}{2+c})^s](1-\beta)+2\beta\delta(2-\mu)]}$.

From above inequalities we have

$\sum_{n=2}^{\infty}n\,a_n\leq\frac{4\beta\delta(1-\mu)}{[\left(\frac{(c_1,m)....(c_t,m)}{(m,m)(b_1,m)...(b_r,m)}\right)(\frac{1+c}{2+c})^s]((1-\beta)+2\beta\delta(2-\mu))}=\gamma$,

thus $f\in N_{n,\gamma}(e)$.

□

**Definition 3.2:**the function $f$(z) defined by (1.7) is said to be a member of the subclass $S_m^{s,c}(\mu,\beta,\delta)$. If there exists a function g $\in S_m^{s,c}(\mu,\beta,\delta)$ such that $\left|\frac{f(z)}{g(z)}-1\right|\leq 1-\zeta, z\in U, 0\leq\zeta<1.$

**Theorem 7.**Let g$\in S_m^{s,c}(\mu,\beta,\delta)$ and

$\zeta=1-\frac{\gamma}{2}d(\mu,\beta,\delta)$. (3.4)

Then $N_{n,\gamma}$(g)c $S_m^{s,c}(\mu,\beta,\delta)$ when $0<\beta\leq 1, 0\leq\mu<\frac{1}{2\delta}, \frac{1}{2}<\delta\leq 1, 0\leq\zeta<1$ and

$d(\mu,\beta,\delta)=\frac{[\frac{(c_1,m)....(c_t,m)}{(m,m)(b_1,m)...(b_r,m)})(\frac{1+c}{2+c})^s]((1-\beta)+2\beta\delta(2-\mu))}{[\frac{(c_1,m)....(c_t,m)}{(m,m)(b_1,m)...(b_r,m)}(\frac{1+c}{2+c})^s]((1-\beta)+2\beta\delta(2-\mu))-2\beta\delta(1-\mu)}$.

**Proof.** Let F$\in N_{n,\gamma}$(g). Then by (3.3)we get $\sum_{n=2}^{\infty}n\,|a_n-b_n|\leq\gamma$,

then

$\sum_{n=2}^{\infty}|a_n-b_n|\leq\frac{\gamma}{2}$.

$\sum_{n=2}^{\infty}b_n\leq\frac{2\beta\delta(1-\mu)}{[\frac{(c_1,m)....(c_t,m)}{(m,m)(b_1,m)...(b_r,m)})(\frac{1+c}{2+c})^s]((1-\beta)+2\beta\delta(2-\mu))}$,

therefore ,

$\left|\frac{F(z)}{g(z)}-1\right|\leq\frac{\sum_{n=2}^{\infty}|a_n-b_n|}{1-\sum_{n=2}^{\infty}b_n}$

$\leq\frac{\gamma}{2}\left[\frac{[\frac{(c_1,m)....(c_t,m)}{(m,m)(b_1,m)...(b_r,m)})(\frac{1+c}{2+c})^s]((1-\beta)+2\beta\delta(2-\mu))}{[\frac{(c_1,m)....(c_t,m)}{(m,m)(b_1,m)...(b_r,m)}(\frac{1+c}{2+c})^s]((1-\beta)+2\beta\delta(2-\mu))-2\beta\delta(1-\mu)}\right]$

$=\frac{\gamma}{2}d(\mu,\beta,\delta)=1-\zeta.$

Then by definition 3.2 we have $f\in S_m^{s,c}(\mu,\beta,\delta)$ .□

**Theorem 8:** Let $f(z)$ *and* $g(z)\in S_m^{s,c}(\mu,\beta,\delta)$ be of the form such that

$f(z)=z-\sum_{n=2}^{\infty}a_n z^n$ and g(z) = z - $\sum_{n=2}^{\infty}b_n z^n$ ,

when $a_n\geq 0$, $b_n\geq 0$.

Then , the Hadamrd product h (z) defiend by $h(z) = z - \sum_{n=2}^{\infty} a_n b_n z^n$ is in the sub class $S_m^{s,c}(\mu,\beta,c)$ when

$$\mu_2 \le [((n-1)(1-\beta)+2\beta\delta(n-\mu_1))^2 \left[\frac{(c_1,m)_{n-1}\dots(c_t,m)_{n-1}}{(m,m)_{n-1}(b_1,m)_{n-1}\dots(b_r,m)_{n-1}})(\frac{1+c}{n+c})^s\right] -2\beta\delta(1-\mu_1)^2(n-1)(1-\beta)-(2\beta\delta)^2(1-\mu_1)^2 n]/[((n-1)(1-\beta) +2\beta\delta(n-\mu_1))^2 \frac{(c_1,m)_{n-1}\dots(c_t,m)_{n-1}}{(m,m)_{n-1}(b_1,m)_{n-1}\dots(b_r,m)_{n-1}})(\frac{1+c}{n+c})^s - (2\delta\beta)^2(1-\mu_1)^2].$$

**Proof.** By theorem 1 , we get

$$\sum_{n=2}^{\infty} \frac{\frac{(c_1,m)_{n-1}\dots(c_t,m)_{n-1}}{(m,m)_{n-1}(b_1,m)_{n-1}\dots(b_r,m)_{n-1}})(\frac{1+c}{n+c})^s((n-1)(1-\beta)+2\beta\delta(n-\mu_1))}{2\beta\delta(1-\mu_1)} a_n \le 1. \quad (3.5)$$

And

$$\sum_{n=2}^{\infty} \frac{\frac{(c_1,m)_{n-1}\dots(c_t,m)_{n-1}}{(m,m)_{n-1}(b_1,m)_{n-1}\dots(b_r,m)_{n-1}})(\frac{1+c}{n+c})^s((n-1)(1-\beta)+2\beta\delta(n-\mu_1))}{2\beta\delta(1-\mu_1)} b_n \le 1. \quad (3.6)$$

We get only to find the lagest $\mu_2$ such that.

$$\sum_{n=2}^{\infty} \frac{\frac{(c_1,m)_{n-1}\dots(c_t,m)_{n-1}}{(m,m)_{n-1}(b_1,m)_{n-1}\dots(b_r,m)_{n-1}})(\frac{1+c}{n+c})^s((n-1)(1-\beta)+2\beta\delta(n-\mu_2))}{2\beta\delta(1-\mu_2)} a_n b_n \le 1.$$

Now by Cauchy –Schwarz inequality , we get

$$\sum_{n=2}^{\infty} \frac{\frac{(c_1,m)_{n-1}\dots(c_t,m)_{n-1}}{(m,m)_{n-1}(b_1,m)_{n-1}\dots(b_r,m)_{n-1}})(\frac{1+c}{n+c})^s((n-1)(1-\beta)+2\beta\delta(n-\mu_1))}{2\beta\delta(1-\mu_1)} \sqrt{a_n b_n} \le 1. \quad (3.7)$$

We need only to show that

$$\frac{\frac{(c_1,m)_{n-1}\dots(c_t,m)_{n-1}}{(m,m)_{n-1}(b_1,m)_{n-1}\dots(b_r,m)_{n-1}})(\frac{1+c}{n+c})^s((n-1)(1-\beta)+2\beta\delta(n-\mu_2))}{2\beta\delta(1-\mu_2)} a_n b_n \le \frac{\frac{(c_1,m)_{n-1}\dots(c_t,m)_{n-1}}{(m,m)_{n-1}(b_1,m)_{n-1}\dots(b_r,m)_{n-1}})(\frac{1+c}{n+c})^s((n-1)(1-\beta)+2\beta\delta(n-\mu_1))}{2\beta\delta(1-\mu_1)} \sqrt{a_n b_n}$$

equivalently

$$\sqrt{a_n b_n} \le \frac{2\beta\delta(1-\mu_2)}{\frac{(c_1,m)_{n-1}\dots(c_t,m)_{n-1}}{(m,m)_{n-1}(b_1,m)_{n-1}\dots(b_r,m)_{n-1}})(\frac{1+c}{n+c})^s((n-1)(1-\beta)+2\beta\delta(n-\mu_2))} \times \frac{\frac{(c_1,m)_{n-1}\dots(c_t,m)_{n-1}}{(m,m)_{n-1}(b_1,m)_{n-1}\dots(b_r,m)_{n-1}})(\frac{1+c}{n+c})^s((n-1)(1-\beta)+2\beta\delta(n-\mu_1))}{2\beta\delta(1-\mu_1)}.$$

But from (3.7)we get

$$\sqrt{a_n b_n} \le \frac{2\beta\delta(1-\mu_1)}{\frac{(c_1,m)_{n-1}\dots(c_t,m)_{n-1}}{(m,m)_{n-1}(b_1,m)_{n-1}\dots(b_r,m)_{n-1}})(\frac{1+c}{n+c})^s((n-1)(1-\beta)+2\beta\delta(n-\mu_1))}$$

Consequently , we also need to prove that

$$\frac{2\beta\delta(1-\mu_1)}{\frac{(c_1,m)_{n-1}\dots(c_t,m)_{n-1}}{(m,m)_{n-1}(b_1,m)_{n-1}\dots(b_r,m)_{n-1}})(\frac{1+c}{n+c})^s((n-1)(1-\beta)+2\beta\delta(n-\mu_1))} \le \frac{(((n-1)(1-\beta)+2\beta\delta(n-\mu_1))}{(((n-1)(1-\beta)+2\beta\delta(n-\mu_2)))}$$

Or , equivalently , that

$$\mu_2 \le [-2\beta\delta(1-\mu 1)^2(n-1)(1-\beta)+((n-1)(1-\beta)+2\delta\beta(n-\mu_1))^2 \times \frac{(c_1,m)_{n-1}\dots(c_t,m)_{n-1}}{(m,m)_{n-1}(b_1,m)_{n-1}\dots(b_r,m)_{n-1}})(\frac{1+c}{n+c})^s - (2\beta\delta)^2(1-\mu_1)^2 n]/[((n-1)(1-\beta)+2\beta\delta(n-\mu_1))^2 \times (\frac{(c_1,m)_{n-1}\dots(c_t,m)_{n-1}}{(m,m)_{n-1}(b_1,m)_{n-1}\dots(b_r,m)_{n-1}})(\frac{1+c}{n+c})^s - (2\delta\beta)^2(1-\mu_1)^2].$$

□

**Theorem9.** Let $f(z) \in S_m^{s,c}(\mu,\beta,\delta)$ be define by (1.7) and $q > -1$. Then the function $G(z)$ defined as $G(z) = \frac{q+1}{z^q}\int_0^z \omega^{q-1} f(\omega), q > -1$,also belongs to$S_m^{s,c}(\mu,\beta,\delta)$

**Proof** : By virtue of G (z) it follows from (1 .7)that

$$G(z) = \frac{q+1}{z^q}\int_0^{\omega}(\omega^q - \sum_{n=2}^{\infty} a_n \omega^{n+q-1})d\omega = z - \sum_{n=2}^{\infty}(\frac{q+1}{q+n})\, a_n z^n$$

But

$$\sum_{n=2}^{\infty} \frac{\frac{(c_1,m)_{n-1}\dots(c_t,m)_{n-1}}{(m,m)_{n-1}(b_1,m)_{n-1}\dots(b_r,m)_{n-1}})(\frac{1+c}{n+c})^s((n-1)(1-\beta)+2\beta\delta(n-\mu\ ))}{2\beta\delta(1-\mu\ )} (\frac{q+1}{q+n})\, a_n \le 1.$$

Since $\frac{q+1}{q+n} \le 1$ and by theorem 1 , so the proof is complete

□

**Theorem 10.** Let F(z) ∈ $S_m^{s,c}(\mu,\beta,\delta)$ be defined by (1.7) and

$$f_\alpha(z) = (1-\alpha)z + \alpha\int_0^z \frac{f(\omega)}{\omega} d\omega \quad (\alpha \geqslant 0, z \in U).$$

Then $F_\alpha(z)$is also in $S_m^{s,c}(\mu,\beta,\delta)$ if $0 \le \alpha \le 2$.

**Proof**. Let $f$ defined by (1.7). Then

$$F_\alpha(z) = (1-\alpha)z + \int_0^z \left(\frac{\omega - \sum_{n=2}^{\infty} a_n\omega^n}{\omega}\right) d\omega$$

$$= z - \sum_{n=2}^{\infty} \frac{\alpha a_n z^n}{n}.$$

By theorem 1 and since $(\frac{\alpha}{2} \le 1)$ we get

$$\sum_{n=2}^{\infty} \frac{\frac{(c_1,m)_{n-1}\dots(c_t,m)_{n-1}}{(m,m)_{n-1}(b_1,m)_{n-1}\dots(b_r,m)_{n-1}})(\frac{1+c}{n+c})^s((n-1)(1-\beta)+2\beta\delta(n-\mu\ ))}{2\delta\beta(1-\mu\ )} (\frac{\delta}{n})\, a_n \le \sum_{n=2}^{\infty} \frac{\frac{(c_1,m)_{n-1}\dots(c_t,m)_{n-1}}{(m,m)_{n-1}(b_1,m)_{n-1}\dots(b_r,m)_{n-1}})(\frac{1+c}{n+c})^s((n-1)(1-\beta)+2\beta\delta(n-\mu\ ))}{2\delta\beta(1-\mu\ )} \left(\frac{\delta}{2}\right) a_n \le 1.$$

Then $f_\alpha(Z)$is in $S_m^{s,c}(\mu,\beta,\delta)$.

□

# فئات معينة من دوال أحادية التكافؤ مع معاملات سالبة معرفة بواسطة العامل الخطي

مازن شاكر محمود ، عبد الرحمن سلمان جمعة ، رحيم احمد منصور
*قسم الرياضيات ، كلية التربية للعلوم الصرفة ، جامعة تكريت ، تكريت ، العراق*

**الملخص**
في هذا البحث تم دراسة الفئات الفرعية من الدوال الاحادية التكافؤ مع معاملات سالبة والتي هي معرفة بواسطة العامل الخطي العام الذي قد قدم . وتم الحصول على النتائج في مقدرات المعاملات والتشوهات وضرب هادمارد ونتائج اخرى تم دراستها وكذلك في هذا البحث تم استخدام العديد من الدوال الكلاسيكية الهندسية العليا بحيث تستطيع ان تكون من $S_m^{s,c}(\mu,\beta,\delta)$ معرفةبالعامل الخطي $\mathcal{H}_m^{s,c}(c_i,b_j)$. الهايبر جميترك من نوع م وتتضمن تقديم فئات فرعية.